# Entropy quotients and correct digits in number-theoretic expansions

## Wieb Bosma[1], Karma Dajani [2] and Cor Kraaikamp [3]

*Radboud University Nijmegen, Utrecht University, and University of Technology Delft*

**Abstract:** Expansions that furnish increasingly good approximations to real numbers are usually related to dynamical systems. Although comparing dynamical systems seems difficult in general, Lochs was able in 1964 to relate the relative speed of approximation of decimal and regular continued fraction expansions (almost everywhere) to the quotient of the entropies of their dynamical systems. He used detailed knowledge of the continued fraction operator. In 2001, a generalization of Lochs' result was given by Dajani and Fieldsteel in [7], describing the rate at which the digits of one number-theoretic expansion determine those of another. Their proofs are based on covering arguments and not on the dynamics of specific maps. In this paper we give a dynamical proof for certain classes of transformations, and we describe explicitly the distribution of the number of digits determined when comparing two expansions in integer bases. Finally, using this generalization of Lochs' result, we estimate the unknown entropy of certain number theoretic expansions by comparing the speed of convergence with that of an expansion with known entropy.

## 1. Introduction

The goal of this paper is to compare the number of digits determined in one expansion of a real number as a function of the number of digits given in some other expansion. For the moment the *regular continued fraction* (RCF) expansion and the *decimal* expansion will serve as examples, the digits then being *partial quotients* and *decimal digits*. It may seem difficult at first sight to compare the decimal expansion with the continued fraction expansion, since the dynamics of these expansions are very different. However, in the early nineteen-sixties G. Lochs [11] obtained a surprising and beautiful result that will serve as a prototype for the results we are after. In it, the number $m(n)$ of partial quotients (or continued fraction digits) determined by the first $n$ decimal digits is compared with $n$. Thus knowledge of the first $n$ decimal digits of the real number $x$ determines completely the value of exactly $m(n)$ partial quotients; to know unambiguously more partial quotients requires more decimal digits. Of course $m(n) = m(x, n)$ will depend on $x$. But asymptotically the following holds (in this paper 'almost all' statements are with respect to the Lebesgue measure $\lambda$).

**Theorem 1.1 (Lochs).** *For almost all $x$:*

$$\lim_{n \to \infty} \frac{m(n)}{n} = \frac{6 \log 2 \log 10}{\pi^2} = 0.97027014 \cdots .$$

[1]Mathematisch Instituut, Radboud Universiteit Nijmegen, Postbus 9010, 6500 GL Nijmegen, The Netherlands, e-mail: bosma@math.ru.nl, url: http://www.math.ru.nl/~bosma/

[2]Fac. Wiskunde en Informatica, Universiteit Utrecht, Postbus 80.000 3508 TA Utrecht, The Netherlands, e-mail: dajani@math.uu.nl, url: http://www.math.uu.nl/people/dajani/

[3]EWI (DIAM), Technische Universiteit Delft, Mekelweg 4, 2628 CD Delft, The Netherlands, e-mail: c.kraaikamp@ewi.tudelft.nl, url: http://ssor.twi.tudelft.nl/~cork/cor.htm

Roughly speaking, this theorem tells us that usually around 97 partial quotients are determined by 100 decimal digits. By way of example Lochs calculated (on an early main-frame computer [12]) that the first 1000 decimals of $\pi$ determine 968 partial quotients.

Examining the function $m(n)$ more closely, one can formulate Lochs' result as follows. Let $B_n(x)$ be the decimal cylinder of order $n$ containing $x$, and $C_m(x)$ the continued fraction cylinder of order $m$ containing $x$, then $m(n)$ is the largest index for which $B_n(x) \subseteq C_m(x)$. Denoting by $S$ the decimal map $S(x) = 10x \mod 1$, and by by $T$ the continued fraction map $T(x) = \frac{1}{x} - \lfloor \frac{1}{x} \rfloor$, then Lochs' theorem says that

$$\lim_{n \to \infty} \frac{m(x)}{n} = \frac{h(S)}{h(T)} \quad \text{(a.e.)}.$$

In 2001, the second author and Fieldsteel generalized in [7] Lochs' result to any two sequences of interval partitions satisfying the result of Shannon–McMillan–Breiman theorem, (see [3], p. 129), that can be formulated as follows.

**Theorem 1.2 (Shannon–McMillan–Breiman).** *Let $T$ be an ergodic measure preserving transformation on a probability space $(X, \mathcal{B}, \mu)$ and let $P$ be a finite or countably infinite generating partition for $T$ for which $H_\mu(P) < \infty$, where $H_\mu(P)$ denotes the entropy of the partition $P$. Then for $\mu$-almost every $x$:*

$$\lim_{n \to \infty} \frac{-\log \mu(P_n(x))}{n} = h_\mu(T).$$

*Here $h_\mu(T)$ denotes the entropy of $T$ and $P_n(x)$ denotes the element of the partition $\bigvee_{i=0}^{n-1} T^{-i}P$ containing $x$.*

The proof in [7] is based on general measure-theoretic covering arguments, and not on the dynamics of specific maps. However, with the technique in [7] it is quite difficult to have a grip on the distribution of $m(n)$. For that, one needs the dynamics of the underlying transformations as reflected in the way the partitions are refined under iterations.

In this paper we give a dynamical proof for certain classes of transformations, and we describe explicitly the distribution of $m(n)$ when comparing two expansions in integer bases. We end the paper with some numerical experiments estimating unknown entropies of certain number theoretical expansions by comparing them with expansions of known entropy.

## 2. Fibred systems

In order to generalize Lochs' result we start with some notations and definitions that place maps like the continued fraction transformation in a more general framework; see also [15], from which the following is taken with a slight modification.

**Definitions 2.1.** Let $B$ be an interval in $\mathbb{R}$, and let $T : B \to B$ be a surjective map. We call the pair $(B, T)$ a *fibred system* if the following conditions are satisfied:

(a) there is a finite or countably infinite set $D$ (called the *digit set*), and
(b) there is a map $k : B \to D$ such that the sets

$$B(i) = k^{-1}\{i\} = \{x \in B : k(x) = i\}$$

form a partition of $B$; we assume that the sets $B(i)$ are intervals, and that $\mathcal{P} = \{B(i) : i \in D\}$ is a generating partition. Moreover, we require that

(c) the restriction of $T$ to any $B(i)$ is an injective continuous map.

The *cylinder of rank $n$ determined by the digits $k_1, k_2, \ldots, k_n \in D$* is the set

$$B(k_1, k_2, \ldots, k_n) = B(k_1) \cap T^{-1}B(k_2) \cap \ldots \cap T^{-n+1}B(k_n);$$

by definition $B$ is the cylinder of rank 0. Two cylinder sets $B_n, B_n^*$ of rank $n$ are called *adjacent* if they are contained in the same cylinder set $B_{n-1}$ of rank $n-1$ and their closures have non-empty intersection:

$$\overline{B_n} \cap \overline{B_n^*} \neq \emptyset.$$

For any fibred system $(B, T)$ the map $T$ can be viewed as a *shift map* in the following way. Consider the following correspondence

$$\Psi : x \mapsto (k_1, k_2, \ldots, k_n, \ldots),$$

where $k_i = j$ if and only if $T^{i-1}x \in B(j)$, then

$$T : \Psi(x) \mapsto (k_2, k_3, \ldots, k_n, \ldots).$$

We write $B_n(x)$ for the cylinder of rank $n$ which contains $x$, i.e.,

$$B_n(x) = B(k_1, k_2, \ldots, k_n) \quad \Longleftrightarrow \quad \Psi(x) = (k_1, k_2, \ldots, k_n, \ldots).$$

A sequence $(k_1, \ldots, k_n) \in D^n$ is called *admissible* if there exists $x \in B$ such that $\Psi(x) = (k_1, \ldots, k_n, \ldots)$.

It will be clear that $\Psi$ associates to any $x \in B$ the 'expansion' $(k_1, k_2, \ldots)$, and that the expansion is determined by $T$ and the digit map $k$. In order to arrive at the desired type of result for the fibred system $(B, T)$, we need to impose some restrictions on $T$.

**Definition 2.2.** A fibred system $(B, T)$ is called a *number theoretic* fibred system and the map $T$ a *number theoretic fibred map* if $T$ satisfies the additional conditions:

(d) $T$ has an invariant probability measure $\mu_T$ which is equivalent to $\lambda$, i.e., there exist $c_1$ and $c_2$, with $0 < c_1 < c_2 < \infty$, such that

$$c_1\lambda(E) \leq \mu_T(E) \leq c_2\lambda(E),$$

for every Borel set $E \subset B$, and

(e) $T$ is ergodic with respect to $\mu_T$ (and hence ergodic with respect to the Lebesgue measure $\lambda$).

**Examples 2.3.** The following examples are used throughout the text.

(i) The first example is the *regular continued fraction*. Here $B = [0, 1) \subset \mathbb{R}$, the map $T$ is given by $T(x) = \frac{1}{x} - k(x)$, and $k(x) = \lfloor \frac{1}{x} \rfloor$. Hence $B(i) = (\frac{1}{i+1}, \frac{1}{i}]$, for all $i \geq 1$. It is a standard result that the map is ergodic with the Gauss-measure as invariant measure. The entropy of $T$ equals $h(T) = \pi^2/(6 \log 2) = 2.373 \cdots$.

(ii) The second standard example is the *$g$-adic expansion*, for an integer $g \geq 2$. Here the map is $T_g(x) = g \cdot x - k(x)$, and $k(x) = \lfloor g \cdot x \rfloor$. Hence $T_g(x)$ is the representative of $g \cdot x \mod 1$ in $B = [0, 1)$, and $B(i) = \lfloor \frac{i}{g}, \frac{i+1}{g} )$, for $0 \leq i < g$. Again this is a number theoretic fibred map, with invariant measure $\lambda$. The entropy of $T_g$ equals $h(T_g) = \log g$.

(iii) A third example is given by the Lüroth expansion. Now $B(i) = [\frac{1}{i+1}, \frac{1}{i})$ for $i \geq 1$. If $T$ is taken as an increasing linear map onto $[0,1)$ on each $B(i)$, the *Lüroth series* expansion is obtained; if $T$ is taken to be decreasing linear onto $[0,1)$ it generates the *alternating Lüroth series* expansion of a number in $[0,1)$. See [1] and [10], where it was shown that the invariant measure is the Lebesgue measure and that the entropy equals

$$h(T) = \sum_{k=1}^{\infty} \frac{\log(k(k+1))}{k(k+1)} = 2.046 \cdots.$$

(iv) As a fourth example we take Bolyai's expansion; again $B = [0,1)$, and now $T : [0,1) \to [0,1)$ is defined by

$$T(x) = (x+1)^2 - 1 - \varepsilon_1(x),$$

where

$$\varepsilon_1(x) = \begin{cases} 0, & \text{if } x \in B(0) = [0, \sqrt{2}-1); \\ 1, & \text{if } x \in B(1) = [\sqrt{2}-1, \sqrt{3}-1); \\ 2, & \text{if } x \in B(2) = [\sqrt{3}-1, 1). \end{cases}$$

Setting $\varepsilon_n = \varepsilon_n(x) = \varepsilon_1(T^{n-1}(x))$, $n \geq 1$, one has

$$x = -1 + \sqrt{\varepsilon_1 + \sqrt{\varepsilon_2 + \sqrt{\varepsilon_3 + \cdots}}}.$$

It is shown in [14] that $T$ is a number theoretic fibred map. Neither the invariant measure $\mu$ nor the exact value of the entropy $h(T)$ is known; compare Experiment (6.3) and the reference to [9] below.

(v) As a final example, let $\beta > 1$ satisfy $\beta^2 = \beta + 1$, i.e., $\beta$ is the *golden mean*. We describe the *$\beta$-continued fraction expansion* of $x$. It is an interesting combination of regular continued fractions and the $\beta$-adic expansion. A *$\beta$-integer* is a real number of the form

$$a_n \beta^n + a_{n-1} \beta^{n-1} + \cdots + a_0$$

where $a_n, a_{n-1}, \ldots, a_0$ is a finite sequence of 0's and 1's without consecutive 1's, i.e., $a_i \cdot a_{i-1} \neq 1$ for $i = 1, 2, \ldots, n$. In fact, $\beta$-integers can be defined for any $\beta > 1$; see e.g. [6] or [8].

For $x \in (0,1)$, let $\lfloor \frac{1}{x} \rfloor_\beta$ denote the largest $\beta$-integer $\leq \frac{1}{x}$, and consider the transformation $T : [0,1) \to [0,1)$, defined by $T(x) = \frac{1}{x} - \lfloor \frac{1}{x} \rfloor_\beta$ for $x \neq 0$, and $T(0) = 0$. Iteration of $T$ generates continued fraction expansions of the form

$$x = \cfrac{1}{b_1 + \cfrac{1}{b_2 + \cfrac{1}{b_3 + \ldots}}},$$

where the $b_i$ are $\beta$-integers. This expansion is called the $\beta$-continued fraction of $x$.

Clearly, this continued fraction expansion is not an $f$-expansion (cf. [15]) in the classical sense, but it has many of the essential properties of $f$-expansions, and thus Rényi's approach in [14] suggests that $T$ has a finite invariant measure equivalent to the Lebesgue measure, with a density bounded away from zero and infinity. Very little seems to be known about such continued fractions. Recently, Bernat showed in [2] that if $x \in \mathbb{Q}(\beta)$, then $x$ has a finite $\beta$-expansion.

We have the following lemma.

**Lemma 2.4.** *Let $T$ be a number theoretic fibred map on $B$, then for almost all $x$:*

$$\lim_{n \to \infty} \frac{\log \lambda(B_n(x))}{\log \mu_T(B_n(x))} = 1.$$

*Proof.* Since

$$c_1 \lambda(B_n(x)) \leq \mu_T(B_n(x)) \leq c_2 \lambda(B_n(x)),$$

it follows that

$$\frac{-\log c_2 + \log \mu_T(B_n(x))}{\log \mu_T(B_n(x))} \leq \frac{\log \lambda(B_n(x))}{\log \mu_T(B_n(x))} \leq \frac{-\log c_1 + \log \mu_T(B_n(x))}{\log \mu_T(B_n(x))}.$$

Taking limits the desired result follows, since $\mathcal{P}$ is a generating partition and therefore $\lim_{n \to \infty} \mu_T(B_n(x)) = 0$ almost surely. □

The following formalizes a property that will turn out to be useful in comparing number theoretic maps.

**Definition 2.5.** Let $T$ be a number theoretic fibred map, and let $I \subset B$ be an interval. Let $m = m(I) \geq 0$ be the largest integer for which we can find an admissible sequence of digits $a_1, a_2, \ldots, a_m$ such that $I \subset B(a_1, a_2, \ldots, a_m)$. We say that $T$ is *r-regular* for $r \in \mathbb{N}$, if for some constant $L \geq 1$ the following hold:

(i) for every pair $B, B^*$ of adjacent cylinders of rank $n$:

$$\frac{1}{L} \leq \frac{\lambda(\overline{B_n})}{\lambda(\overline{B_n^*})} \leq L,$$

(ii) for almost every $x \in I$ there exists a positive integer $j \leq r$ such that either $B_{m+j}(x) \subset I$ or $B_{m+j}^* \subset I$ for an adjacent cylinder $B_{m+j}^*$ of $B_{m+j}(x)$.

Thus $r$-regularity expresses, loosely speaking, that if the expansion $T$ agrees to the first $m$ digits in both endpoints of a given interval $I$, then for a.e $x \in I$ there exists a cylinder of rank at most $m + r$ with the property that it or an adjacent cylinder contains $x$, and is contained entirely in $I$; moreover, the size of two adjacent cylinders differs by not more than the constant factor $L$.

**Examples 2.6.** We give three examples.

(i) Let $I \subset [0, 1)$ be a subinterval of positive length, and let $m = m(I)$ be such that $B_m$ is the smallest RCF-cylinder containing $I$, i.e., there exists a vector $(a_1, \ldots, a_m) \in \mathbb{N}^m$ such that $B_m = B_m(a_1, \ldots, a_m) \supset I$ and $\lambda(I \cap ([0, 1) \setminus B_{m+1}(a_1, \ldots, a_m, a))) > 0$ for all $a \in \mathbb{N}$. One can easily check that $r = 3$ in case of the RCF.

Furthermore, it is well-known, see e.g. (4.10) in [3], p. 43, that

$$\lambda(B_m(a_1, \ldots, a_m)) = \frac{1}{Q_m(Q_m + Q_{m-1})},$$

which yields, together with the well-known recurrence relations for the partial fraction denominator sequence $(Q_m)_{m\geq 0}$ (see [3], (4.2), p. 41):

$$Q_{-1} := 0; \quad Q_0 := 1; \quad Q_m = a_m Q_{m-1} + Q_{m-2}, \quad m \geq 1,$$

that

$$\lambda(B_m(a_1,\ldots,a_m)) \leq 3\lambda(B_m(a_1,\ldots,a_m+1)),$$

for all $(a_1,\ldots,a_m) \in \mathbb{N}^m$, i.e., the continued fraction map $T$ is 3-regular, with $L = 3$.

(ii) The $g$-adic expansion is not $r$-regular for any $r$. Although adjacent rank $n$ cylinders are of the same size, so $L = 1$ can be taken in (2.5)(i), property (2.5)(ii) does not hold. This can be seen by taking $x$ in a very small interval $I$ around $1/g$.

(iii) The *alternating* Lüroth map is 3-regular with $L = 2$.

## 3. A comparison result

We would like to compare two expansions. In this section we will therefore assume that $S$ and $T$ are both number theoretic fibred maps on $B$. We denote the cylinders of rank $n$ of $S$ by $A_n$, and those of $T$ by $B_n$. We let $m(x,n)$ for $x \in B$ and $n \geq 1$ be the number of $T$-digits determined by $n$ digits with respect to $S$, so $m(x,n)$ is the largest positive integer $m$ such that $A_n(x) \subset B_m(x)$. Yet another way of putting this, is that both endpoints of $A_n(x)$ agree to exactly the first $m$ digits with respect to $T$.

We have the following general theorem.

**Theorem 3.1.** *Suppose that $T$ is $r$-regular. Then for almost all $x \in B$*

$$\lim_{n\to\infty} \frac{m(x,n)}{n} = \frac{h(S)}{h(T)}.$$

*Proof.* For $x \in B$ let the first $n$ $S$-digits be given; then by definition of $m = m(x,n)$ one has that $A_n(x) \subset B_m(x)$. Since $T$ is $r$-regular, there exists $1 \leq j \leq r$ such that

$$B_{m+j} \subset A_n(x) \subset B_m(x),$$

where $B_{m+j}$ is either $B_{m+j}(x)$, or adjacent to it of the same rank. By regularity again, for some $L \geq 1$,

$$\frac{1}{L}\lambda(B_{m+r}(x)) \leq \frac{1}{L}\lambda(B_{m+j}(x)) \leq \lambda(B_{m+j}),$$

so

$$\frac{1}{n}\left(-\log L + \log \lambda(B_{m+r}(x))\right) \leq \frac{1}{n}\log \lambda(A_n(x)) \leq \frac{1}{n}\log \lambda(B_m(x)),$$

and the result follows from Lemma (2.4) (applied to both $T$ and $S$), and the Theorem of Shannon-McMillan-Breiman (1.2). □

**Example 3.2.** Let $m_g^{\text{RCF}}(x,n)$ be the number of partial quotients of $x$ determined by the first $n$ digits of $x$ in its $g$-adic expansion. Then for almost all $x$:

$$\lim_{n\to\infty} \frac{m_g^{\text{RCF}}(x,n)}{n} = \frac{6\log 2\log g}{\pi^2}.$$

This generalizes Lochs' theorem to arbitrary $g$-adic expansions.

**Example 3.3.** Let $m_g^{\text{ALE}}(x, n)$ be the number of alternating Lüroth digits of $x$ determined by the first $n$ digits of $x$ in its $g$-adic expansion. Then for almost all $x$:

$$\lim_{n \to \infty} \frac{m_g^{\text{ALE}}(x, n)}{n} = \frac{\log g}{\sum_{k=1}^{\infty} \frac{\log k(k+1)}{k(k+1)}}.$$

Note that the conditions in Theorem 3.1 are not symmetric in $S$ and $T$. Since the regularity condition is not satisfied by $h$-adic expansions, Theorem 3.1 is of no use in comparing $g$-adic and $h$-adic expansions, nor in proving that the first $n$ regular partial quotients determine usually $\pi^2/(6 \log 2 \log g)$ digits in base $g$.

## 4. Comparing radix expansions

In this section we compare $g$-adic and $h$-adic expansions, for integers $g, h \geq 2$, by explicitly studying the distribution of $m(n)$.

Given a positive integer $n$, there exists a unique positive integer $\ell = \ell(n)$, such that

$$h^{-(\ell+1)} \leq g^{-n} \leq h^{-\ell}. \tag{4.1}$$

Thus the measure $\lambda(A_n)$ of a $g$-cylinder of rank $n$ is comparable to that of an $h$-cylinder of rank $\ell(n)$.

It follows that

$$\limsup_{n \to \infty} \frac{\ell(n)}{n} \leq \frac{\log g}{\log h} \leq \liminf_{n \to \infty} \frac{\ell(n) + 1}{n},$$

i.e.,

$$\lim_{n \to \infty} \frac{\ell(n)}{n} = \frac{\log g}{\log h}, \tag{4.2}$$

which is the ratio of the entropies of the maps $T_g$ and $T_h$ introduced in (2.3)(ii). Of course, one expects the following to hold:

$$\lim_{n \to \infty} \frac{m_g^{(h)}(x, n)}{n} = \frac{\log g}{\log h} \quad \text{(a.e.)}. \tag{4.3}$$

which is the analog of Lochs' result for the maps $T_g$ and $T_h$; here $m_g^{(h)}(x, n)$ is defined as before: it is the largest positive integer $m$ such that $A_n(x)$ is contained in the $h$-adic cylinder $B_m(x)$. One should realize that, in general, $\ell(n)$ has no obvious relation to $m(n)$; in fact equation (4.2) is merely a statement about the relative speed with which $g$-adic and $h$-adic cylinders shrink.

Let $x \in [0, 1)$ be a generic number for $S = T_g$ for which we are given the first $n$ digits $t_1, t_2, \ldots, t_n$ of its $g$-adic expansion. These digits define the $g$-adic cylinder

$$A_n(x) = A_n(t_1, t_2, \ldots, t_n).$$

Let $m(n) = m_g^{(h)}(n)$ and $\ell(n)$ be as defined above; note that

$$m(n) \leq \ell(n), \quad \text{for all } n \geq 1. \tag{4.4}$$

The sequence $(m(n))_{n=1}^{\infty}$ is non-decreasing, but may remain constant some time; this means that it 'hangs' for a while, so $m(n+t) = \cdots = m(n+1) = m(n)$, after which it 'jumps' to a larger value, so $m(n+t+1) > m(n+t)$. Let $(n_k)_{k \geq 1}$ be the subsequence of $n$ for which $m(n)$ 'jumps', that is, for which $m(n_k) > m(n_k - 1)$.

**Lemma 4.5.** *For almost all* $x$,

$$\lim_{k \to \infty} \frac{m_g^{(h)}(n_k)}{n_k} = \frac{\log g}{\log h}.$$

*Proof.* By definition, $B_{m(n_k)}(x) \supset A_{n_k}(x)$. The cylinder $B_{m(n_k)}(x)$ consists of $h$ cylinders of rank $m(n_k) + 1$, and $A_{n_k}(x)$ intersects (at least) two of these, since otherwise some $B_{m(n_k)+1} \supset A_{n_k}(x)$, contradicting maximality of $m(n_k)$, and thus it contains an endpoint $e$ of some $B_{m(n_k)+1}$ lying in the interior of $B_{m(n_k)}(x)$.

On the other hand, by definition of $n_k$, we know that $m(n_k) > m(n_k - 1)$, so $B_{m(n_k)}(x) \not\supset A_{n_k-1}(x)$, and therefore $A_{n_k-1}(x)$ contains an endpoint $f$ of $B_{m(n_k)}(x)$ as well as $e$. Now $e$ and $f$ are at least $\lambda(B_{m(n_k)+1}) = h^{-(m(n_k)+1)}$ apart. Therefore

$$h^{-(m(n_k)+1)} \leq \lambda(A_{n_k-1}(x)) = g \cdot \lambda(A_{n_k}(x)) = g \cdot g^{-n_k} \leq g \cdot h^{-\ell(n_k)},$$

by (4.1). Hence

$$h^{-(m(n_k)+1)} \leq g \cdot h^{-\ell(n_k)},$$

which, in combination with (4.4) implies

$$\ell(n_k) - 1 - \frac{\log g}{\log h} \leq m(n_k) \leq \ell(n_k), \quad k \geq 1. \tag{4.6}$$

But from (4.6) the Lemma follows immediately. $\square$

Next we would like to show that (4.3) follows from (4.6); this is easy in case $h = 2$. Compare [7].

**Corollary 4.7.** *With notations as above, for any* $g \in \mathbb{N}$, $g \geq 2$

$$\lim_{n \to \infty} \frac{m_g^{(2)}(x, n)}{n} = \frac{\log g}{\log 2} \quad \text{(a.e.)}.$$

*Proof.* If $h = 2$ the mid-point $\xi$ of $B_{m(n)}$ lies somewhere in $A_n$. Now $m(n+1) = m(n)$ if this mid-point $\xi$ is located in $A_{n+1}(t_1, t_2, \ldots, t_n, t_{n+1})$, that is, if it is located in the same $g$-adic cylinder of order $n + 1$ as $x$. Notice that this happens with probability $\frac{1}{g}$, and that the randomness here is determined by the $g$-adic digit $t_{n+1}$. To be more precise, let $H$ denote the event that we will 'hang' at time $n$, i.e., the event that $m(n) = m(n+1)$, and let $D$ be a random variable with realizations $t \in \{0, 1, \ldots, g-1\}$, defined by $\xi \in A_{n+1}(t_1, \ldots, t_n, D)$. Now

$$\mathrm{P}(H) = \sum_{i=0}^{g-1} \mathrm{P}(H|D = i) \cdot \mathrm{P}(D = i),$$

and from $\mathrm{P}(H|D = i) = 1/g$ for $0 \leq i \leq g-1$ it then follows that $\mathrm{P}(H) = 1/g$; due to the discrete uniform distribution of the digit $t_{n+1}$ of $x$ we do not have to know the probabilities $\mathrm{P}(D = i)$.

Thus we see that for each $k \geq 1$ (and with $n_0 = 0$) the random variable

$$v_k = n_k - n_{k-1}$$

is geometrically distributed with parameter $p = 1/g$. Furthermore, the $v_k$'s are independent, and, since $n_k \geq k$

$$1 \leq \frac{n_{k+1}}{n_k} = \frac{n_k + n_{k+1} - n_k}{n_k} = 1 + \frac{1}{n_k} v_{k+1} \leq 1 + \frac{1}{k} v_{k+1}.$$

Because $v_1, v_2, \ldots$ are independent and identically distributed with finite expectation, it follows from the Lemma of Borel-Cantelli that $\lim_{k \to \infty} \frac{v_{k+1}}{k} = 0$ (a.e.), and therefore

$$\lim_{k \to \infty} \frac{n_{k+1}}{n_k} = 1 \quad \text{(a.e.)}.$$

Given any $n \geq 1$, there exists $k = k(n) \geq 0$ such that $n_k < n \leq n_{k+1}$. Since $m(n_k) \leq m(n) \leq m(n_{k+1})$ one has

$$\frac{m(n_{k(n)})}{n_{k(n)+1}} \leq \frac{m(n)}{n} \leq \frac{m(n_{k(n)+1})}{n_{k(n)}},$$

But then the result follows from (4.6) with $h = 2$. $\qquad\qquad\Box$

If $h \geq 3$ the situation is more complicated; there might be more than one midpoint $\xi_i$ 'hitting' $A_n$. In case only one $\xi_i$ 'hits' $A_n$ (as is always the case when $h = 2$), we speak of a TYPE 1 situation. In case there is more than one 'hit' we speak of a TYPE 2 situation. Now change the sequence of jump-times $(n_k)_{k \geq 1}$ by adding those $n$ for which we are in a TYPE 2 situation, and remove all $n_k$'s for which we are in a TYPE 1 situation. Denote this sequence by $(n_k^*)_{k \geq 1}$. If this sequence is finite, we were originally in a TYPE 1 situation for $n$ sufficiently large. In case $(n_k^*)_{k \geq 1}$ is an infinite sequence, notice that for every $n$ for which there exists a $k$ such that $n = n_k^*$ one has 1-regularity, i.e.,

$$B_{m(n_k^*)+1} \subset A_{n_k^*}(x) \subset B_{m(n_k^*)}(x),$$

where $B_{m(n_k^*)+1}$ is either $B_{m(n_k^*)+1}(x)$ or an adjacent interval. Following the proof of Theorem 3.1 one has

$$\lim_{k \to \infty} \frac{m(n_k^*)}{n_k^*} = \frac{\log g}{\log h} \quad \text{(a.e.)}.$$

Let $(\hat{n}_k)_{k \geq 1}$ be the sequence we get by merging $(n_k)_{k \geq 1}$ and $(n_k^*)_{k \geq 1}$. Notice that

$$\lim_{k \to \infty} \frac{m(\hat{n}_k)}{\hat{n}_k} = \frac{\log g}{\log h} \quad \text{(a.e.)}. \tag{4.8}$$

If $(\tilde{n}_k)_{k \geq 1}$ is $\mathbb{N} \setminus (\hat{n}_k)_{k \geq 1}$, we are left to show that

$$\lim_{k \to \infty} \frac{m(\tilde{n}_k)}{\tilde{n}_k} = \frac{\log g}{\log h} \quad \text{(a.e.)}.$$

Let $n = \tilde{n}_j$ for some $j \geq 1$. Then there exist unique $k = k(j)$ and $h = h(j)$ such that $n_{k+1} = \hat{n}_{h+1}$, $n_k \leq \hat{n}_h$ and

$$\hat{n}_h < n < n_{k+1}.$$

Since $n_{k+1} - \hat{n}_h = \hat{n}_{h+1} - \hat{n}_h = v_{h+1}$ is geometrically distributed with parameter $1/g$, and since $v_1, v_2, \ldots$ are independent with finite expectation, it follows as in the case $h = 2$, that

$$\frac{\hat{n}_{h+1}}{\hat{n}_h} = 1 + \frac{1}{\hat{n}_h} v_{h+1} \leq 1 + \frac{1}{h} v_{h+1} \to 1 \quad \text{as } j \to \infty \quad \text{(a.e.)}.$$

But then (4.3) follows from (4.8) and from

$$\frac{m(\hat{n}_h)}{\hat{n}_{h+1}} \leq \frac{m(n)}{n} \leq \frac{m(\hat{n}_{h+1})}{\hat{n}_h}.$$

We have proved the following theorem; again, see also [7].

**Theorem 4.9.** *Let* $g, h \in \mathbb{N}_{\geq 2}$. *Then for almost all* $x$

$$\lim_{n \to \infty} \frac{m_g^{(h)}(x,n)}{n} = \frac{\log g}{\log h}.$$

## 5. Remarks on generalizations

At first sight one might think that the proof of Theorem 4.9 can easily be extended to more general number theoretic fibered maps, by understanding the distribution of $m(n)$. However, closer examination of the proof of Theorem 4.9 shows that there are two points that make generalizations of the proof hard, if not impossible. The first is the observation, that the 'hanging time' $v_k$ (in a TYPE 1 situation) is geometrically distributed, and that the $v_k$'s are independently identically distributed. Recall that this observation follows from the fact that the digits given by the $g$-adic map $S_g$ are independently identically distributed and have a discrete uniform distribution. As soon as this last property no longer holds (e.g., if $S = T_\gamma$, the expansion with respect to some non-integer $\gamma > 1$), we need to know $\mathrm{P}(D = i)$ for each $i \in \{0, 1, \ldots, \lfloor \gamma \rfloor - 1\}$, and this might be difficult.

Even if we assume $S$ to be the $g$-adic map, another problem arises when $T$ is not the $h$-adic map, but, for example, $T = T_\gamma$. In that case all the ingredients of the proof of Theorem 4.9 seem to work with the exception that at the 'jump-times' $(n_k)_k$ we might not be able to show that

$$\ell(n_k) - C \leq m(n_k) \leq \ell(n_k), \tag{5.1}$$

with $C \leq 1$ some fixed constant. Notice that from (5.1) it would follow that

$$\lim_{k \to \infty} \frac{m(n_k)}{n_k} = \frac{\log g}{\log \gamma}.$$

For one class of $\gamma \in (1, 2)$ one can show that (5.1) still holds. These $\gamma$'s are the so-called 'pseudo-golden mean' numbers; $\gamma > 1$ is a 'pseudo-golden mean' number if $\gamma$ is the positive root of $X^k - X^{k-1} - \cdots - X - 1 = 0$, for some $k \in \mathbb{N}, k \geq 2$ (in case $k = 2$ one has that $\gamma = \beta = \frac{1}{2}(1 + \sqrt{5})$, which is the golden mean). These 'pseudo-golden mean' numbers $\gamma$ are all Pisot numbers, and satisfy

$$1 = \frac{1}{\gamma} + \frac{1}{\gamma^2} + \cdots + \frac{1}{\gamma^k}.$$

**Theorem 5.2.** *With notations as before, let* $S = T_g$ *and* $T = T_\gamma$, *with* $g \in \mathbb{N}_{\geq 2}$, *and* $\gamma > 1$ *a 'pseudo golden number'. Then*

$$\lim_{n \to \infty} \frac{m_g^{(\gamma)}(n)}{n} = \frac{\log g}{\log \gamma}, \quad \text{(a.e.)}.$$

This result is immediate from [7]. Here we present a different proof, based on the distribution of the number of correct digits, $m(n)$.

*Proof of Theorem 5.2.* In case $k = 2$, a digit 1 is always followed by a zero, and only the $\gamma$-cylinders $B_m^{(\gamma)}$ corresponding to a sequence of $\gamma$-digits ending with 0 are refined. Thus from the definition of $m(n)$ it follows that the last digit of $B_{m(n)}^{(\gamma)}(x)$ is always 0 (if it were 1, the choice was wrong). Note that in this case

$$\lambda(B_{m(n)}^{(\gamma)}(x)) = \gamma^{-m(n)}.$$

Let $\ell(n)$ be defined as before, and notice that at a 'jump-time' $n = n_k$ one has

$$\gamma^{-(m(n)+2)} \leq g \cdot g^{-n} \leq g \cdot \gamma^{-\ell(n)}, \tag{5.3}$$

from which

$$\ell(n) - 2 - \frac{\log g}{\log \gamma} \leq m(n).$$

Since $m(n) \leq \ell(n)$ we see that (5.1) follows in case $\gamma$ equals the golden mean.

In case $k \geq 3$ the situation is slightly more complicated; let us consider here $k = 3$. Now any sequence of $\gamma$-digits ending with two consecutive 1's must be followed by a zero, and therefore — by the definition of $m(n)$ — the last digit of $B_{m(n)}^{(\gamma)}(x)$ is either 0, or is 1 which is preceded by 0. One can easily convince oneself then that

$$\lambda(B_{m(n)}^{(\gamma)}(x)) = \begin{cases} \gamma^{-m(n)} & \text{if } d_{m(n)}(x) = 0, \\ \gamma^{-(m(n)+1)} + \gamma^{-(m(n)+2)} & \text{if } d_{m(n)}(x) = 1. \end{cases}$$

Here $d_{m(n)}(x)$ is the $m(n)^{\text{th}}$ $\gamma$-digit of $x$.

Now let $n = n_k$ be a 'jump-time'. If $d_{m(n)}(x) = 0$, then $B_{m(n)}^{(\gamma)}(x)$ consists of two $\gamma$-cylinders, one of length $\gamma^{-(m(n)+1)}$ and one of length $\gamma^{-(m(n)+2)} + \gamma^{-(m(n)+3)}$. One of these two cylinders is contained in $A_{n-1}(x)$ (due to the fact that $n$ is a 'jump-time'), and (5.3) is therefore satisfied. Of course one has $m(n) \leq \ell(n)$ in this case. In case $d_{m(n)}(x) = 1$ we see that $B_{m(n)}^{(\gamma)}(x)$ consists of two $\gamma$-cylinders, one of length $\gamma^{-(m(n)+1)}$ and one of length $\gamma^{-(m(n)+2)}$. Now

$$\lambda(B_{m(n)}^{(\gamma)}(x)) = \gamma^{-(m(n)+1)} + \gamma^{-(m(n)+2)} \leq \gamma^{-m(n)},$$

and by definition of $\ell(n)$ one has $m(n) \leq \ell(n)$. Again one of these two sub-cylinders of $B_{m(n)}^{(\gamma)}(x)$ is contained in $A_{n-1}(x)$, and one has

$$\gamma^{-(m(n)+2)} \leq g \cdot g^{-n} \leq g \cdot \gamma^{-\ell(n)}.$$

We see that (5.3) is again satisfied. For $k \geq 4$ the proof is similar; one only needs to consider more cases.                                                                    □

## 6. Estimation of (unknown) entropies

In this section we report on some numerical experiments. Our experiments were carried out using the computer algebra system MAGMA, see [4]. The general set-up of the experiments was as follows. We choose $n$ random digits for the expansion of a number $x$ in $[0, 1)$ with respect to a number theoretic fibred map $S$, and compute $m_S^T(x, n)$, the number of digits of $x$ with respect to the number theoretic fibred map $T$, determined completely by the first $n$ digits of $x$ with respect to $S$. This is done by comparing the $T$-expansions of both endpoints of the $S$-cylinder $A_n(x)$.

**Experiment 6.1.** The first experiment was designed to test the set-up. It aimed to verify Lochs' result: we let $S = T_{10}$, the decimal expansion map, and $T = T^{\text{RCF}}$ for the regular continued fraction map. For $N = 1000$ random real numbers of $n = 1000$ decimal digits, we computed the number of partial quotients determined by the endpoints of the decimal cylinder $A_n(x)$.

We ran this experiment twice; the averages (over 1000 runs) were 970.534 and 969.178 (with standard deviations around 24). Compare these averages to the value $970.270 \cdots$ predicted by Theorem 1.1.

**Experiment 6.2.** In this experiment we compared the relative speed of approximation of $g$-adic and $h$-adic expansions:

| $h \to$ | 2 | | | 7 | | | 10 | | |
|---|---|---|---|---|---|---|---|---|---|
| $g \downarrow$ | pred | observed | | pred | observed | | pred | observed | |
| 2 | - | - | - | **0.356** | 0.355 | 0.355 | **0.301** | 0.300 | 0.300 |
| 7 | **2.807** | 2.805 | 2.805 | - | - | - | **0.845** | 0.844 | 0.844 |
| 10 | **3.322** | 3.320 | 3.320 | **1.183** | 1.182 | 1.182 | - | - | - |

In fact we compared binary, 7-adic and decimal expansions of $N = 1000$ random real numbers with each of the other two expansions, again by computing to how many $h$-digits both endpoints of the $g$-adic cylinder $A_{1000}(x)$ agreed. The two tables list the values found in each of two rounds of this experiment and compares the result with the value found by Theorem 4.9. Note that by (4.4) always $m(n) \leq \ell(n)$ and thus the observed ratio should approximate $\log g / \log h$ from below. By reversing the roles of $g$ and $h$, it is thus possible in this case to approximate the entropy quotient by sandwiching!

**Experiment 6.3.** By choosing random real numbers of 1000 decimal digits, we attempted to estimate the unknown entropy of the Bolyai map (cf. (2.3)(iv)). We found that on average 1000 digits determined 2178.3 Bolyai digits in the first run of $N = 250$ random reals, and 2178.0 in the second; the standard deviations were around 13. These values would correspond to an entropy of around 1.0570 or 1.0572.

Choosing $N$ random Bolyai expansions of length 1000 (by finding the decimal expansion to 1000 digits for random reals), we determined the number of decimal digits in one experiment, and the number of regular continued fraction partial quotients in another experiment. In the first run of $N = 1000$, we found that on average 458 decimal digits were determined, in the second run 457.8 on average (with standard deviations around 4). These values would correspond to an entropy of around 1.0545 for the first, and 1.0541 for the second run, according to Theorem 5 in [7]. In runs with $N = 1000$ of the second experiment we found that 444.36 and 444.37 partial quotients were determined on average, with standard deviations of around 17. This suggests an entropy of 1.0545 or 1.0546 according to Theorem 3.1.

These entropy estimates for the Bolyai-map were reported in a preliminary version of this paper [5]. Using a different method, based on the fixed-points of the Bolyai-map, Jenkinson and Pollicott were able to find in [9] sharper estimates for the entropy. To be precise, they found 1.056313074, and showed that this is within $1.3 \times 10^{-7}$ of the true entropy.

**Experiment 6.4.** In our final experiment we attempted to estimate the entropy of the $\beta$-continued fraction map, see (2.3)(v).

As before, we chose real numbers with 1000 random decimal digits and determined the number of $\beta$-integer partial quotients determined by them. In the first run, the 1000 decimals determined 877.922 partial quotients on average, in the second run 878.125 partial quotients were determined. These values correspond to an entropy of 2.021 and 2.022.